\documentclass{svjour2}                    % onecolumn

\smartqed  % flush right qed marks, e.g. at end of proof
 \usepackage{mathptmx}
\usepackage{pdfsync}
\usepackage{amsmath,amssymb,bbm}
\usepackage{graphicx,graphics}
\usepackage{color}
\newcommand{\ts}{\hspace{0.5pt}}

\newcommand{\ee}{{\rm e}}
\newcommand{\dd}{\,\mathrm{d}}

\newcommand{\NN}{\mathbb{N}}
\newcommand{\EE}{\mathbb{E}}
\newcommand{\PP}{\mathbb{P}}
\newcommand{\one}{\mathbbm{1}}

\DeclareMathOperator{\Geo}{Geo}

\begin{document}

\title{Lines of descent under selection\thanks{This article is dedicated to the memory of Hans-Otto Georgii,
whose joint work with the first author on ancestral lines in multitype branching processes \cite{Georgii_Baake_03,Baake_Georgii_07} laid foundations and provided
motivation for the line of research reviewed here.}
}


\author{Ellen Baake         \and
        Anton Wakolbinger %etc.
}


\institute{E. Baake \at
              Faculty of Technology, Bielefeld University,  \\
              Tel.: +49-521-106-4896
              \email{ebaake@techfak.uni-bielefeld.de}           %  \\
           \and
           A. Wakolbinger \at
            Institute of Mathematics, Goethe-Universit\"at Frankfurt am Main, \\
            Tel.: +49-69-798-28651
            \email{wakolbinger@math.uni-frankfurt.de} % \\
}

\date{Received: date / Accepted: date}

\maketitle

\begin{abstract}
We review recent progress on ancestral processes related to mutation-selection models, both in the deterministic and the stochastic setting.
We  mainly rely on two concepts, namely, the \emph{killed ancestral selection graph} and the \emph{pruned lookdown ancestral selection graph}. The killed ancestral selection graph gives a representation of the type of a random individual from a stationary population, based upon the individual's potential ancestry back until the mutations that define the individual's type. The pruned lookdown ancestral selection graph allows one to trace  the ancestry of individuals from a stationary distribution
back into the \emph{distant} past, thus leading to the stationary distribution of \emph{ancestral} types. We illustrate the results by applying them  to a prototype model for the error threshold phenomenon.

%Insert your abstract here. Include keywords, PACS and mathematical
%subject classification numbers as needed.
\keywords{mutation-selection model \and killed ancestral selection graph \and pruned lookdown ancestral selection graph \and error threshold}
\subclass{
60J27 % Continuous-time Markov processes on discrete state spaces
\and 60J75 %jump processes
\and  92D15 %Problems related to evolution
\and  05C80 %random graphs
}
\end{abstract}

\section{Introduction}
\label{intro}
Understanding the interplay of   mutation and
selection is a major topic of population genetics research. Among statistical physicists,
the deterministic Crow-Kimura model \cite{Crow_Kimura_56,Crow_Kimura_70} is particularly well known: It describes the  \emph{parallel} action of
mutation and selection on the genetic composition of  an effectively
infinite population. 
%Here independence means that mutation events are not coupled to reproduction, but may happen any time during the life cycle of an organism. 
The population is identified with a probability distribution on some type space,  and the dynamics is  defined via a system of ordinary differential equations (ODEs).
A particularly relevant type space is $\{0,1\}^\ell$, the set of all binary sequences of  length $\ell$. Starting in the 1980s, Leuth\"ausser \cite{Leuthaeusser_86,Leuthaeusser_87} and Tarazona \cite{Tarazona_92} established a connection between a discrete-time version of the quasispecies model \cite{Eigen_71} (a
mutation-selection model  on sequence space where mutation happens on the occasion of reproduction, independently at every site of the sequence with probability $p$) and
a classical Ising model. More precisely, the mutation-reproduction matrix governing the  dynamical system was identified with the transfer matrix of a two-dimensional Ising model with nearest-neighbour interaction between the rows, but arbitrary interaction within the rows; see \cite{Baake_Gabriel_00} for a review.  This connection made mutation-selection models very popular in the statistical physics community.
Interest mainly focussed on the stationary type distribution, that is, on the balance between mutation and selection.
In the late 1990s, a connection was established between the  Crow-Kimura sequence space model (where sites mutate independently at the same rate $\mu$) and an Ising quantum chain \cite{Baake_Baake_Wagner_97}. Here, the matrix governing the ODE system is equivalent to the Hamiltonian of an Ising quantum chain in a transverse field, with arbitrary interactions within the chain; see \cite{Baake_Gabriel_00} for the details.  This fact opened the
toolbox of quantum statistical mechanics for evolution models.   Quantum-mechanical probabilities and expectations played a pivotal role in the calculation of the leading eigenvalue of the mutation-reproduction matrix and hence of the stationary distribution of types; however, these quantum-mechanical objects  seemed to be lacking a
biological interpretation, a fact that remained mysterious for some time. Later, it turned out that the desired translation into classical probabilities emerges if one traces the ancestry of individuals from a stationary distribution back into the distant past, thus leading to the stationary distribution of  \emph{ancestral types}
\cite{Hermisson_etal_02}; see \cite[Appendix~A]{Hermisson_etal_02} for the translation from the quantum-mechanical into the classical picture.
In probability theory, the concept of ancestral type distributions had previously been introduced for multitype branching processes
 \cite{Jagers_Nerman_84,Jagers_89,Jagers_92}, which are closely related to mutation-selection models.  Furthermore, a variational principle was established, which allows one to characterise both the present and the ancestral populations \cite{Baake_Georgii_07,Georgii_Baake_03,Hermisson_etal_02}, and  has been widely used in the physics community, see, e.g., \cite{Baake_Baake_Bovier_Klein_05,Garske_06,Garske_Grimm_04,Leibler_Kussell_10,Peliti_02,Sughiyama_Kobayashi_17}. By now, the equilibrium structure of  deterministic mutation-selection models, both in the forward and backward directions of
time, is quite well understood.

A parallel development, mainly within probability theory and biomathematics, took place in the context of  stochastic processes. This was initiated by the seminal work of Fisher \cite{Fisher_30} and Wright \cite{Wright_31} in the 1930s; further landmarks were the contributions of Mal\'{e}cot in 1948 \cite{Malecot_48}, Feller in 1951 \cite{Feller_51}, and Moran in 1958 \cite{Moran_58}. They   also described the fluctuations due to random reproduction over long time scales, which are absent in the deterministic dynamics; this development is traced in \cite{Nagylaki_89}, \cite[Ch.~1.4.3]{Ewens_04}, and \cite{Baake_Wakolbinger_15}. In the sequel,  the introduction of the coalescent process by  Kingman  in the early 1980s \cite{Kingman_82a,Kingman_82b}  stimulated a wealth of results on genealogies of samples of individuals taken from populations at present. Consequently, the focus of population genetics research  changed from the prospective to the  retrospective point of view. Indeed, understanding  the ancestral processes contributes decisively to understanding the genetic structure of today's populations. While coalescent theory was restricted to the neutral case (that is, the case without selection) for the first 15 years, genealogical constructions that allowed to deal with selection became available later; namely, the ancestral selection graph by Krone and Neuhauser \cite{Krone_Neuhauser_97} and the lookdown construction by Donnelly and Kurtz \cite{Donnelly_Kurtz_99}.

  It is the purpose of this article to review recent progress in this direction, and to bring together some of the lines of research mentioned above. We will start from the stochastic model with mutation and selection, trace back  ancestral lines, and thus obtain insight into both the present and the ancestral distributions of types. In the deterministic limit, we will recover some previously-known results and give them additional meaning. The stochastic models  are yet more challenging than the  deterministic ones  and, so far, the ancestral distributions are only well explored for the case with two types, a beneficial and a deleterious one (encoded by 0 and 1, respectively); the multitype case is the subject of current research. For the sake of a unifying description of the deterministic and stochastic models, we will restrict ourselves to the two-type case throughout.  The paper is organised as follows. We first introduce the model in both its deterministic and stochastic versions, along with some basic facts (Section~\ref{sec:model}); we then lay out the general concept of the ancestral selection graph (ASG, Section~\ref{sec:ASG}). The remainder of the paper is devoted to two recent constructions based on the ASG, namely, the killed ASG (Section~4) and the pruned lookdown ASG (Section~5). The killed ASG allows one to determine the stationary type distribution  in terms of a genealogical picture, both in the deterministic and the stochastic models; the pruned lookdown ASG serves the analogous purpose for the stationary distribution of the ancestral types.  

\section{Model and basic facts}
\label{sec:model}
A widely-used prototype model of population genetics is the
Moran model with two types under mutation and  selection. It assumes
a  population of fixed size $N \in \mathbb{N}$ in which each individual is characterised
by a type $i \in \{ 0,1 \}$. An individual of type $i$ may, at any instant in continuous time, 
do either of two things: it may reproduce (at rate $r^N_i >0$),  or it may mutate (at rate $u^N > 0$ unless stated otherwise); the dependence on $N$ will become important later. As to reproduction, let $r^N_0=1/2+s^N$ with $s^N \geqslant 0$,  and $r^N_1=1/2$; so type $0$ has \emph{selective advantage} $s^N$ and is hence understood as the beneficial type, whereas type $1$ is selectively inferior.
(Since the birth rates of the two types differ, one speaks of \emph{fecundity} or \emph{fertility selection}.)
When an individual reproduces, its single offspring inherits the 
parent's type and replaces a uniformly chosen individual, possibly its own parent.  When an individual mutates, the new  type is $j$ with probability $\nu_j^{}$ with $0 \leqslant \nu_j^{} \leqslant 1$ and $\nu_0^{}
+\nu_1^{} =1$. Note that this includes the possibility of silent mutations, where the type is the same before and after the event.

\begin{figure}[ht]
\begin{center}
%\psfrag{t}{\LARGE$t$}
%\psfrag{0}{\Large$0$}
%\psfrag{1}{\Large$1$}
\includegraphics[angle=0, width=0.56\textwidth]{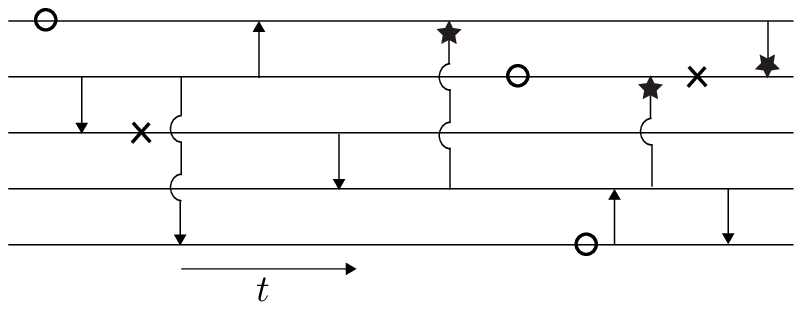} \\[2mm]
\includegraphics[angle=0, width=0.6 \textwidth]{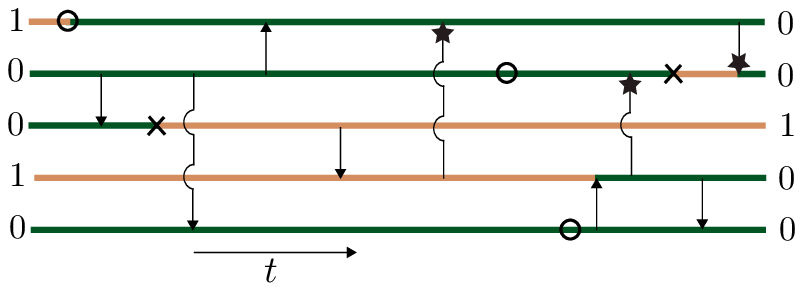}
\caption{The Moran model with two-way mutation and selection.  Crosses represent mutations to type $1$, circles mutations to type $0$. Neutral and selective reproduction events are depicted as arrows with `usual' and star-shaped arrowheads, respectively. Top: untyped Moran particle system; bottom: the same realisation but with types assigned. Dark green: type~0; light brown: type~1.
\label{MoranModel}}
\end{center}
\end{figure}

The Moran model has a well-known graphical illustration as an interacting particle system, as illustrated in Fig.~\ref{MoranModel}.
The individuals are represented by horizontal line pieces, with forward (physical) time $t$ running from left to right in the figure (and all further pictures that show functions of time).
An  arrow indicates a
reproduction event with the parent at its tail and the offspring at its head. 
For later use, we decompose reproduction events into neutral and selective ones. Neutral arrows (with the usual arrowheads in Fig.~\ref{MoranModel}) appear 
at rate $1/2$ per individual and hence at rate $1/(2N)$ per ordered pair of lines; selective arrows (those with a star-shaped arrowhead) appear at rate $s^N/N$ per ordered pair of lines, 
irrespective of the line types. The \emph{type-dependent} reproduction rates $r_i^N$ are then obtained by the convention that neutral arrows
are used by all individuals, whereas  selective  arrows are used by type-$0$ individuals and are ignored otherwise.
Mutations to type $0$ (type~1) appear at rate $u^N \nu_0$ ($u^N \nu_1$) on every line and are marked by circles (crosses).

The above description of the particle system is an example of an \emph{untyped} construction in the following sense. The graphical elements (arrows, circles, crosses) are laid out at constant rates, regardless of the types. The benefit of this is twofold. First, a given realisation of the graphical representation can then be used \emph{with any initial configuration of types}. Second, the law for the graphical elements is simple: they appear at constant rates, in the manner of Poisson processes for every line (or pair of lines, respectively), regardless of the types; in particular, the law is the same in the forward and backward directions of time.  Some events will turn out as silent after the assignment of the initial types  (such as a selective arrow encountered by a type-1 individual, or a circle encountered by a type-0 individual); in this way, one can assign the graphical elements to the lines regardless of type, which will turn out  advantageous in the constructions. Indeed, working with untyped constructions will prove as key to the analysis.

The graphical representation may also be read out in an alternative way that gives rise to \emph{viability selection} rather than fecundity selection. Under viability selection, the \emph{death rates} are type-dependent. More precisely, individuals of type $0$ die at rate $d^N_0=1/2$, whereas individuals of type~$1$ die at rate $d^N_1=1/2+s^N$; when an individual dies, it is replaced by an offspring of an individual that is chosen uniformly from the population. In the particle representation, the effect of the neutral arrows remains the same as before; but for selective arrows, the type of the individual \emph{at the tip of the arrow} now is decisive. If this individual is of type $1$, then the individual at the tail of the arrow places offspring via the arrow; if the individual at the tip is of type $0$, the arrow is ignored.

While the parental relationships may differ between fecundity and viability selection, the behaviour is the same at the level of the graphical representation, in the sense that fecundity and viability selection lead to the same typed picture in the lower panel of Fig.~\ref{MoranModel}. In particular, $X_t^{N}$, the 
proportion of type-$0$ individuals at time $t$ in a population of size $N$, coincides under both modes of selection, for any given  realisation of the particle picture. This is because, under both variants,  $X_t^{N}$ increases   by $1/N$ if and only if an  arrow (neutral or selective) points from a type-$0$ individual to a type-$1$ individual; it decreases by $1/N$ if a neutral arrow points from a type-1 individual to a type-0 individual.

One usually studies the model in an $N \to \infty$ limit. The following two limits are by far the most relevant.

\begin{figure}
\begin{center}
\includegraphics[ width=0.48\textwidth]{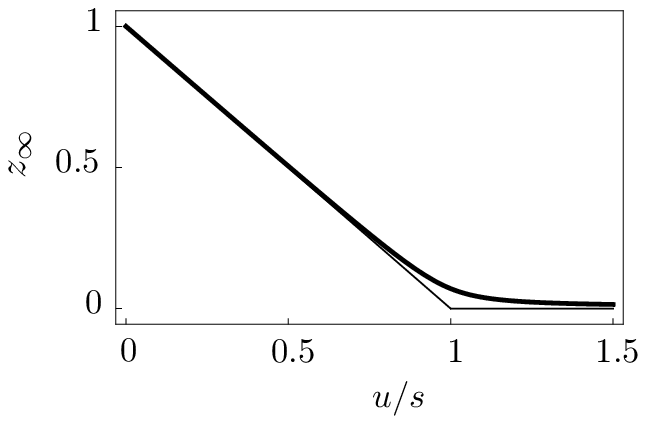} \quad
\includegraphics[ width=0.48\textwidth]{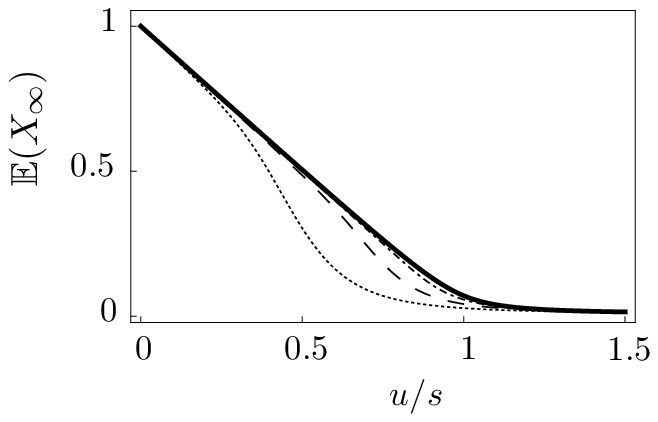} 
\end{center}
\caption{\label{stat_forward}
The stationary proportion of type 0 in the deterministic limit (left) and the corresponding expectation in the diffusion approximation (right) for $s=0.001$ as a function of the mutation rate. Further parameters in left panel: $\nu_0=0$ (thin line),$\nu_0=0.005$ (bold line). In right panel: $\nu_0=0.005$ (for all curves);  $N=10^4$ (dotted), $N=3 \cdot 10^4$ (dashed), $N=10^5$ (dot-dashed); the bold line is again $z_\infty$, as in the left panel.
}
\end{figure}

\paragraph{Law of large numbers (deterministic limit):} 
Here, one lets $N \to \infty$ without any rescaling of parameters and time;
so $s^N \equiv s$, $u^N \equiv u$.
If $ X^N_0 \to x$ as $N \to \infty$,  then
$ \big ( X^N_t \big)_{t \geqslant 0}$ converges weakly to
$\big (z(t) \big )_{t \geqslant 0}$,
where $z(t)$ solves the Riccati differential equation
\begin{equation}\label{det_muse}
 \dot {z}  =  s z (1-z) + u \nu^{}_0 (1-z) - u \nu^{}_1 z
\end{equation}
with initial value $z(0)=x$; 
the solution  is known explicitly, see for instance \cite{Cordero_17b}.
For $0<\nu_0\leqslant 1$, the solution of \eqref{det_muse} converges to 
\begin{equation}\label{z_infty}
z_\infty  := \lim_{t \to \infty} z(t) = \begin{cases} \frac{1}{2} \Big (1-\frac{u}{s} + \sqrt{\Big ( 1-\frac{u}{s}\Big )^2 + 4 \frac{u}{s} \nu^{}_0 } \Big ), & s>0 \\
\nu_0, & s=0, \end{cases}
\end{equation}
independently of $x$, hence $z_\infty$ is a globally stable equilibrium for \eqref{det_muse}.
The case with highly asymmetric mutation ($\nu_0 \ll \nu_1$) is of particular interest.
It is  widely used as a prototype for the sequence-space model with a single-peaked landscape. The latter model assumes that one single type (the `wildtype sequence') reproduces at rate $1/2+s$, whereas all  others
(the $2^{\ell}-1$ `mutants') reproduce at rate\footnote{In the deterministic setting, it is more common to assume a neutral reproduction rate of $1$ rather than $1/2$. We work with the rate $1/2$ here in order to obtain the pair coalescence rate of $1$ that is standard in coalescence theory (see Section~\ref{sec:ASG}). Note that  the deterministic dynamics \eqref{det_muse} is unaffected by the neutral reproduction rate anyway.} 1/2, and mutation corresponds to a random walk on the $\ell$-dimensional
hypercube.  The prototype model emerges in the approximation of  lumping all mutants  into a single type,
see \cite{Eigen_McCaskill_Schuster_89} for a review.
In the limiting case  $\nu_0=0$, which results from the prototype model in the limit $\ell \to \infty$, Eq.~\eqref{z_infty} reduces to
\begin{equation}\label{eq_spl}
z^{}_\infty = \begin{cases} 1 - \frac{u}{s},  & u \leqslant s \\ 0, & u>s. \end{cases}
\end{equation}
More precisely,   for $u\leqslant s$  the solution of \eqref{det_muse} converges to $z^{}_\infty= 1 - u/s$ 
 for all $0 < x \leqslant 1$, whereas for $u> s$ it converges to $z^{}_\infty= 0$ 
 even for all $0 \leqslant x \leqslant 1$. 
Thus 
Eq.~\eqref{eq_spl} means that the beneficial type is lost from the population when the mutation rate
surpasses the selective advantage --- a phenomenon that became prominent under the
name of \emph{error threshold} \cite{Eigen_71,Eigen_McCaskill_Schuster_89}, and may be seen as a phase transition. See Fig.~\ref{stat_forward} for an illustration.

This deterministic limit
is a special case of a general dynamical law of large numbers by
 Kurtz \cite{Kurtz_71}, see also \cite[Thm.~11.2.1]{Ethier_Kurtz_86}; in the present case the convergence extends  to the stationary state \cite{Cordero_17b}.
%  Rewriting
%$y = (y^{}_0, y^{}_1) :=  ( z , 1-z )$, we obtain
%\begin{equation}\label{muse_vector}
%% \dot y^{}_i = y^{}_i \big ( r^{}_i - \langle r,y \rangle \big ) + (y M)_i ,
%% \quad i \in \{0,1\} \,,
%% \end{equation}
%% where $ r = (r^{}_0, r^{}_1) = (1/2+s,  1/2)$ holds the reproduction rates,
%% $\langle .,.\rangle$ denotes the scalar product, $M = \begin{pmatrix} -u \nu_0 & u \nu_0 \\ u \nu_1 & - u \nu_1 \end{pmatrix}$
%% is the mutation generator, and $(y M)_i$ is the $i$-th component of the (row) vector $yM$.
%% Eq.~\eqref{muse_vector} is the  (Crow-Kimura) mutation-selection equation mentioned in the Introduction, for the case of two types. 
%% This formulation is useful since it may be transformed into a \emph{linear} system governed by the matrix  $A:= M + \text{diag}(r)$, which has the usual solution via the matrix exponential and leads back to the solution of \eqref{det_muse} by normalisation. In particular, the  equilibrium of Eq.~\eqref{muse_vector} is the (probability-normalised) Perron--Frobenius left eigenvector of $A$. In a similar spirit, $y$ may be understood as the vector of relative frequencies in a two-type branching process with first-moment generator $A$, when the population has grown to infinite size; this fruitful connection has been exploited in \cite{Hermisson_etal_02,Georgii_Baake_03,Baake_Georgii_07}.
 A comprehensive review of deterministic mutation-selection models is provided in \cite{Buerger_00}.

\paragraph{Diffusion limit:} 
Here one assumes that the selective advantage and the mutation rate both depend on the population size $N$, obeying
%Parameters are rescaled such that 
$\lim_{N \to \infty} N s^N = \sigma$ and 
$\lim_{N \to \infty} N u^N = \vartheta$ with $0 \leqslant \sigma, \vartheta < \infty$, and time is sped up by a 
factor of $N$. In the limit $N \to \infty$,  $(X_{tN}^{N})_{t \in \mathbb{R}}$ then converges (in distribution) to  $X:=(X_t)_{t \in \mathbb{R}}$, the  Wright--Fisher diffusion on $[0,1]$ characterised by
the drift coefficient
$\alpha(x) = \sigma x (1-x) +  \vartheta \nu_0^{} (1-x) -  \vartheta \nu_1^{} x$
and the diffusion coefficient
$\beta(x) = (1/2) x(1-x)$. That is, $(X_t)_{t \in \mathbb{R}}$ follows the stochastic differential equation
\begin{equation}\label{sde}
{\rm d} X_t =  \alpha(X_t) \ts {\rm d} t + \sqrt{\beta(X_t)} \ts {\rm d} W_t, 
\end{equation}
where $(W_t)_{t \geqslant 0}$ is standard  Brownian motion. 
Here the drift coefficient captures the deterministic trend; note that it has the same form as the right-hand side of \eqref{det_muse} but with a different scaling of the parameters. The diffusion term captures the fluctuations due to random (neutral) reproduction in the finite population, which persists in the diffusion limit.

For the limiting case $\vartheta=0$, $X$ is absorbing (in $0$ or $1$). The probability of absorption  in state 1 is given by
\[
h(x) := \PP(X \text{ absorbs in } 1 \mid X_0=x) = \frac{1- \exp(-2 \sigma x)}{1-\exp(-2 \sigma )};
\]
this is a classical result of Mal\'{e}cot \cite{Malecot_48} and Kimura \cite{Kimura_62}. Note that the limit $\sigma=0$ renders  the neutral fixation probability $h(x)=x$.
 For $\vartheta>0$ and $\nu_0=1$ ($\nu_0=0$), the process absorbs in 1 (0) with probability one. For $\vartheta>0$ and $0 <\nu_0<1$, the process has a stationary distribution known as Wright's distribution, which has density 
 \begin{equation}\label{pi}
 \pi(x) = C (1-x)^{2 \vartheta \nu^{}_1-1} x^{2 \vartheta \nu^{}_0-1} \ee^{2 \sigma x}, \quad 0<x<1,
 \end{equation}
 where $C$ is a normalising constant. Comprehensive reviews of diffusion models in population genetics may be found in  \cite[Ch.~4, 5]{Ewens_04} and \cite[Ch.~7, 8]{Durrett_08}.

\section{The ancestral selection graph}
\label{sec:ASG}
\begin{figure}[ht]
\begin{center}
%\psfrag{I}{\Large$I$}
%\psfrag{C}{\Large$C$}
%\psfrag{D}{\Large$D$}
%\psfrag{0}{\Large$0$}
%\psfrag{1}{\Large$1$}
\includegraphics[angle=0,width=.45 \textwidth]{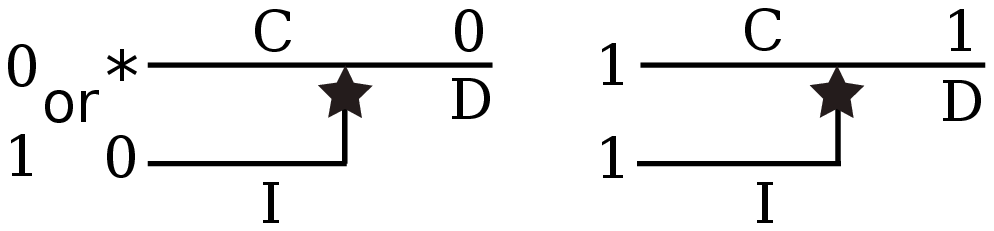} \qquad 
\includegraphics[angle=0,width=.45 \textwidth]{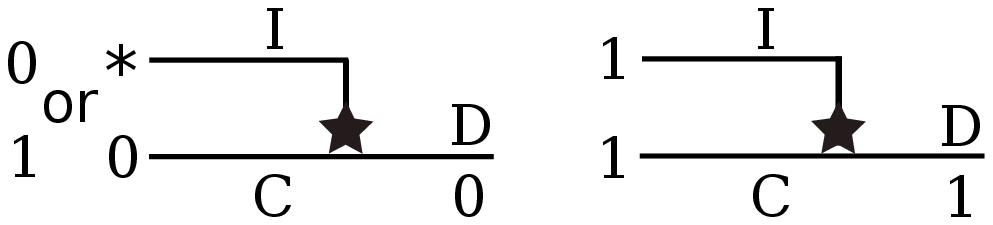}
\caption{The pecking order for fecundity selection (first and second diagram) and viability selection (third and fourth diagram). Incoming branch (I), continuing branch (C), and descendant (D). A $*$ stands for an arbitray  type. The type  combinations for I and C on the left  of each diagram  lead to the type of D noted on the right.  Physical time runs from left to right.
\label{IC}}
\end{center}
\end{figure}

One central concept to study ancestries and  genealogies of (samples of) individuals is the  \emph{ancestral selection graph (ASG)} of Krone and Neuhauser \cite{Krone_Neuhauser_97}. It was originally formulated for the case of fecundity selection; we review it here and include the straightforward extension to viability selection.
A basic principle is to start again with an untyped picture and to
understand selective arrows 
as unresolved reproduction events backward in time. Namely, 
the descendant has two \textit{potential ancestors}, the 
\textit{incoming branch} (at the tail) and the
\textit{continuing branch} (at the tip), see  Fig.~\ref{IC}. 
In the case with fecundity selection, the incoming 
branch is the ancestor if it is of type $0$, otherwise the continuing one is ancestral. With viability selection, the priority 
is interchanged:  The continuing branch is 
the ancestor if it is of type $0$, otherwise the incoming one is ancestral. We will loosely refer to this hierarchy as the \emph{pecking order}.
Note that, in any case, the descendant is of type 1 if and only if both potential parents are of type 1.

The construction of the ASG  starts from a  sample taken at time $t > 0$, to which we refer as the `present'. One first ignores the types and traces back the lines of all \emph{potential ancestors} the individuals may have at times previous to $t$.  In the graphical representation with a finite population size $N$,
a neutral arrow that joins two potential ancestral  lines appears at rate 
$1/N$ per currently extant pair of potential ancestral lines,
then giving rise to a \textit{coalescence event}, i.e.
the two lines merge into a single one, thus reducing the number of potential ancestors by one. In contrast to the neutral case, lines can also \emph{branch} into two. This happens whenever a potential ancestral line
is hit by a selective arrow that
emanates from outside the current set of $n$ potential ancestral lines, that is,  at rate $ s_N^{} (N-n) /N$. Viewed backward in time, the line that is hit  splits into an
incoming and continuing branch as described above. Every extant ancestral line is hit by a selective arrow from
within the current set at rate  $ s_N^{} (n-1)/N$.  Such an event is called a \emph{collision}; it does not change the number of potential ancestors.  Mutation events 
are superposed on the lines of the  ASG
at rates $u^{N} \nu_0^{}$ and $u^{N} \nu_1^{}$, respectively.

From now on, we will concentrate on the $N \to \infty$ limits. In both the diffusion and the deterministic limit,
collision events have probability 0. In the diffusion limit \cite{Krone_Neuhauser_97}, one is left with branching events (at rate $\sigma$ per line),
coalescence events (at rate 1 per  pair of lines), and mutation events (at rate $\vartheta \nu_0^{}$ and $\vartheta \nu_1^{}$ per line, respectively).
In the deterministic limit \cite{Cordero_17}, one  loses the coalescence events, so that only branching (rate $s$ per line) and mutation (rate $u \nu_0^{}$ and $u \nu_1^{}$
per line)
survive. In the diffusion limit, the number of lines in the graph always remains finite (with probability 1), whereas it diverges  in the deterministic limit as time tends to infinity (for any $s>0$). 

\begin{figure}
\includegraphics[ width=0.45\textwidth]{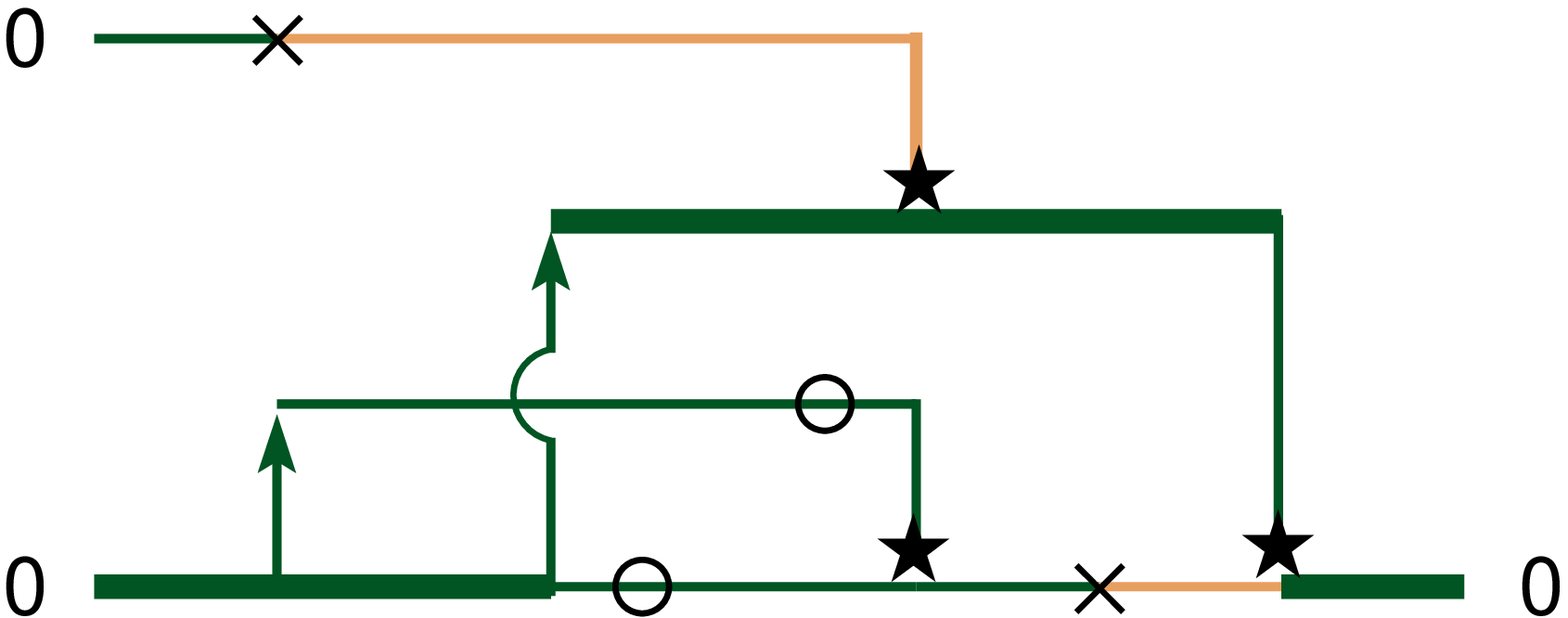} \qquad
\includegraphics[ width=0.45\textwidth]{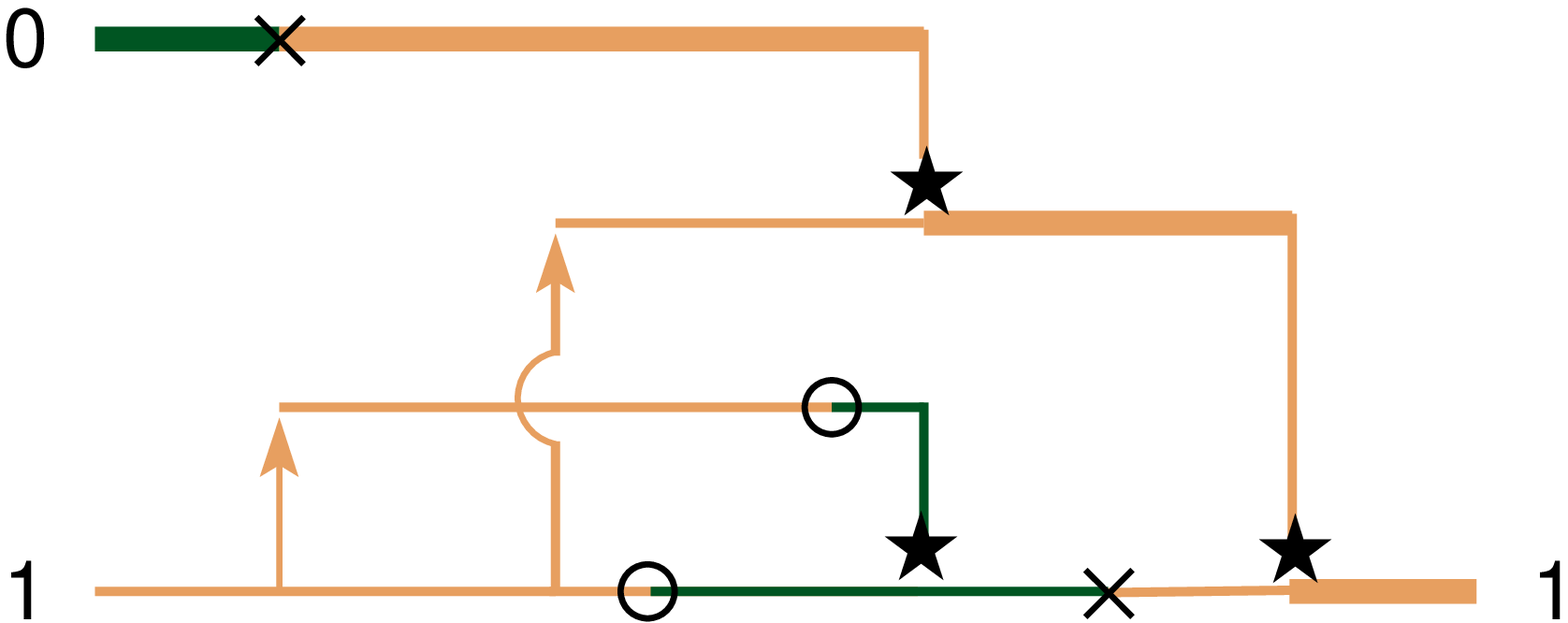} \\[5mm]
\includegraphics[ width=0.45\textwidth]{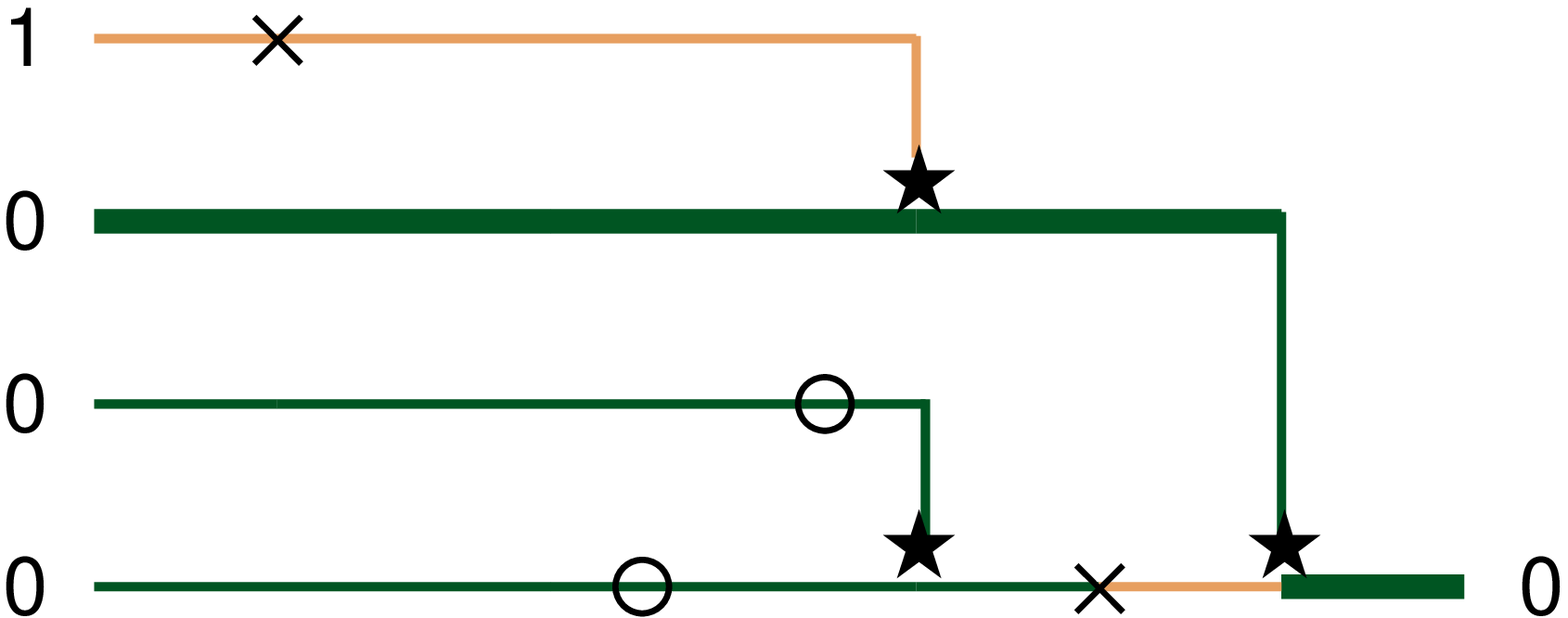} \qquad
\includegraphics[ width=0.45\textwidth]{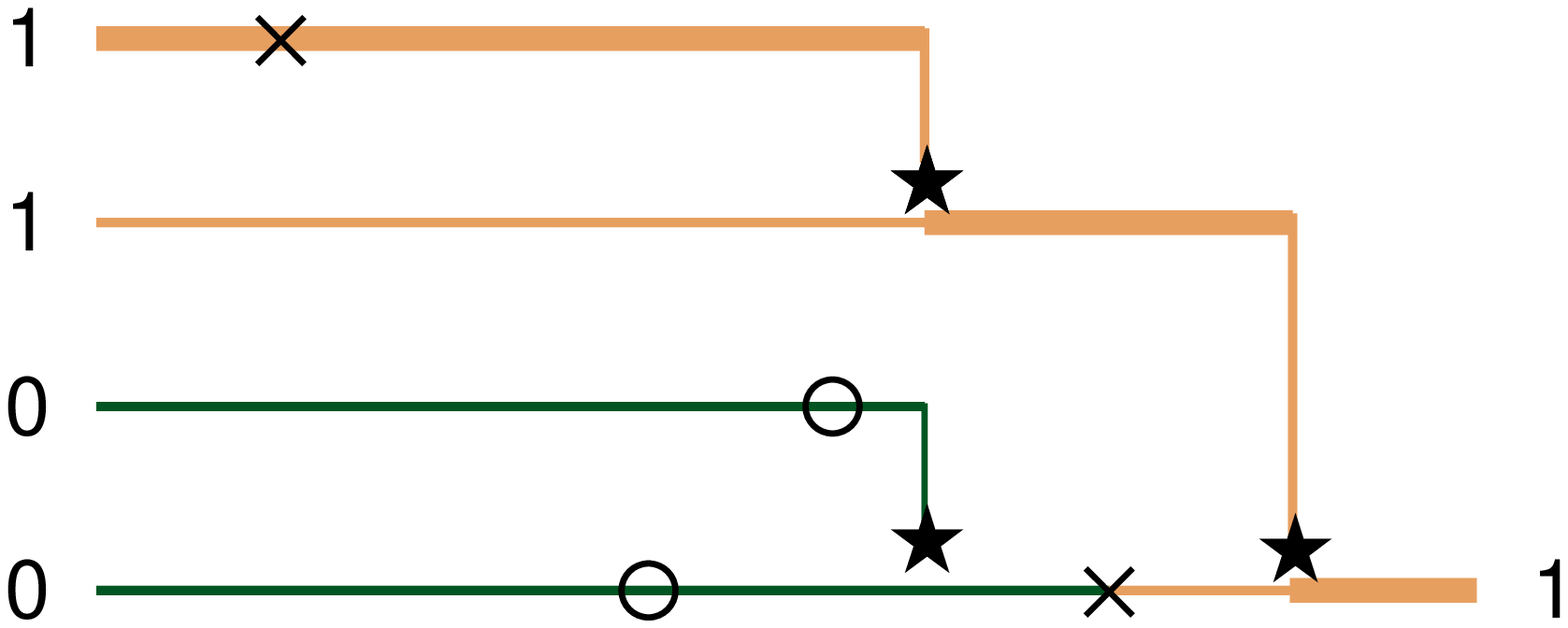} 
\caption{\label{fig:identify_line}  ASG with viability selection, and true ancestral line (bold) for a sample of size 1, in the diffusion limit (top) and the deterministic limit (bottom). The figure shows how different ancestral configurations may lead to different types of the sampled individual; note, however, that different ancestral type configurations may also lead to the same type.  Dark green: type~0; light brown: type~1; circles: beneficial mutations; crosses: deleterious mutations; arrows with `usual' arrowheads: neutral reproduction events; arrows with star-shaped arrowheads: selective reproduction events. 
}
\end{figure}

When the ASG (in either of the two limits) has been constructed backward in time until time 0, say, then one assigns types to its lines  in the  time interval $[0, t]$.
Given the frequency $x$ of the beneficial type at time $0$, one  first draws the types of the lines at time $0$ independently and identically distributed (i.i.d.) according to the probability vector $(x,1-x)$; more precisely, every line is assigned type $0$ with probability $x$ and type $1$ with probability $1-x$, independently of each other, and independently of the ASG.  One then propagates 
the types forward in time, respecting the mutation
events.
In this way, the (backward in time) branching events may now be resolved to reveal who is the true parent in every single case. Removing the non-ancestral branches then gives the true genealogy of the sample.  Fig.~\ref{fig:identify_line} shows some examples for the case of a sample of size 1; in this case, there is always exactly one true ancestral line. Note that we order the lines in a specific way: with viability selection, the incoming branch is always placed immediately above the continuing line in a branching event; with fecundity selection, the incoming branch is placed  just below the continuing branch. In a coalescence event, the ancestral line continues on the lower of the two lines, so that the arrow points upwards. The latter is an element of the lookdown construction \cite{Donnelly_Kurtz_99}. This ordering is allowed since both the initial assignment of types and the  dynamics of the ASG  are invariant under permutation of lines (that is, they are \emph{exchangeable}); it will ease the graphical representation of the true ancestral line later on.

\section{The killed ASG and the stationary type distribution}
\label{sec:killed}
Let us now explain how the ASG can provide insight into the  type distribution at present. To this end, we will introduce the \emph{killed ASG}, separately for the two limits. We need not distinguish between fecundity and viability selection here, since we will only be interested in the type of the descendant at any given branching event; one therefore need not decide whether the incoming or the continuing line is ancestral (as in Fig.~\ref{IC}). Nevertheless, we adhere to the  ordering of the lines described in the previous section.

\begin{figure}
\begin{center}
\includegraphics[ width=0.22\textwidth]{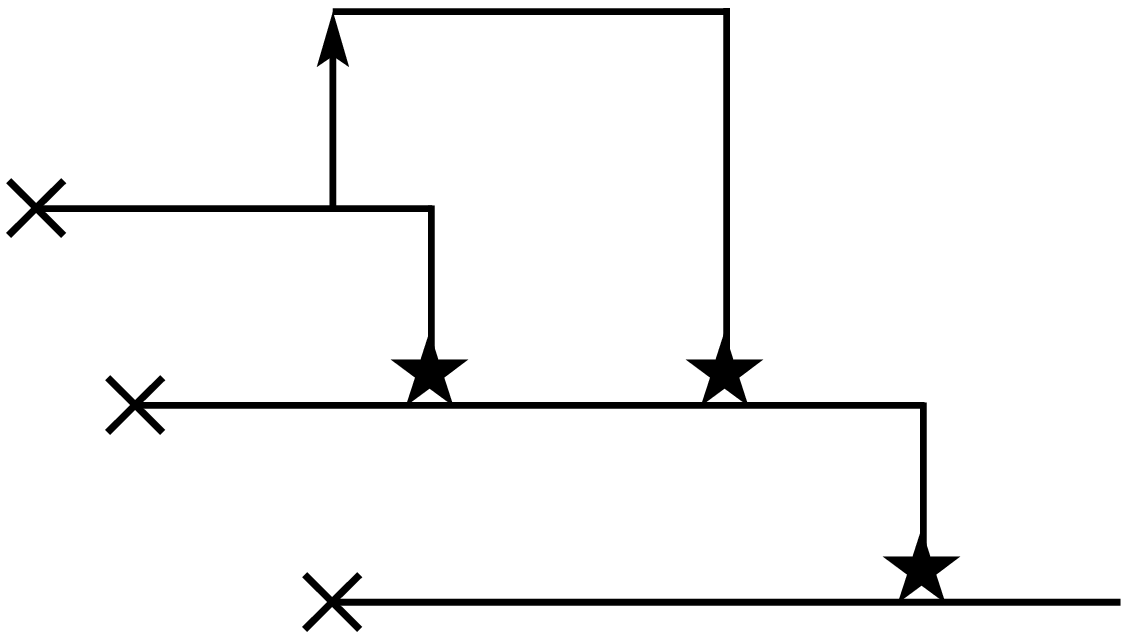} \quad
\includegraphics[ width=0.22\textwidth]{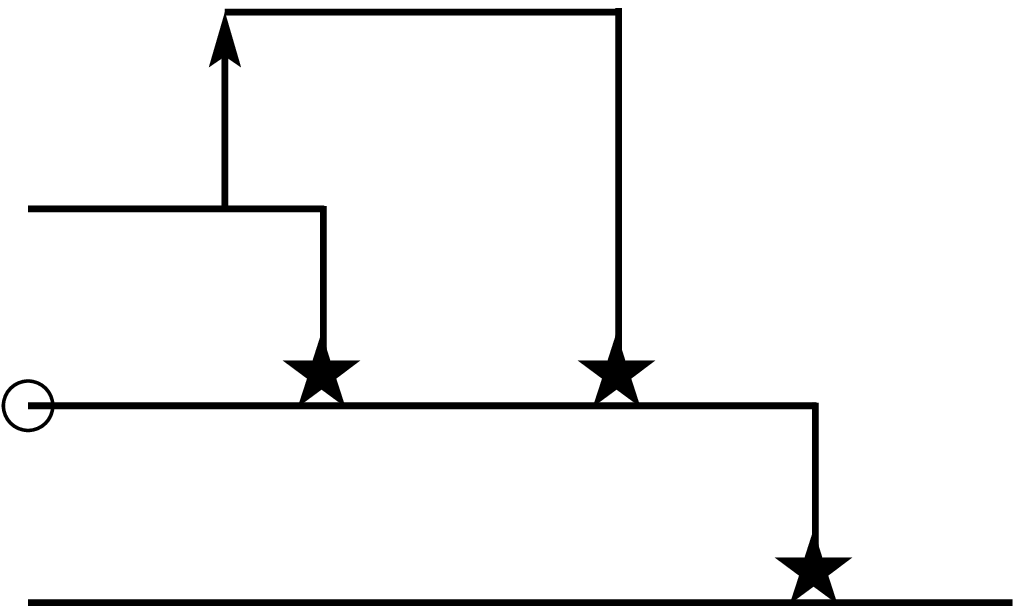} \quad
\includegraphics[ width=0.22\textwidth]{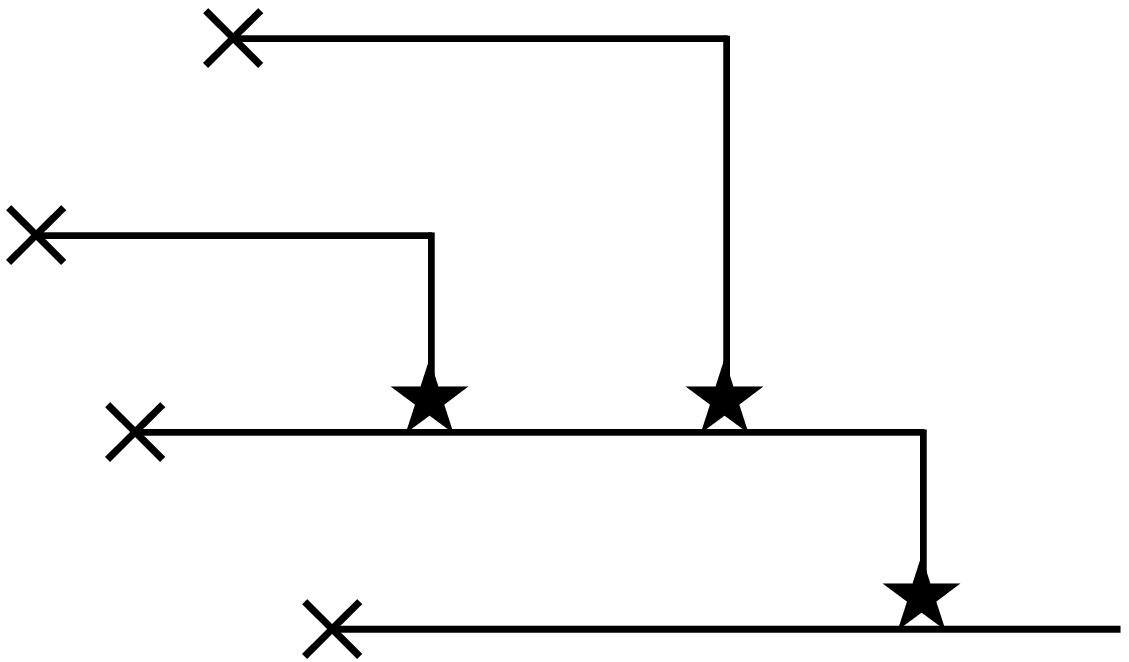} \quad
\includegraphics[ width=0.22\textwidth]{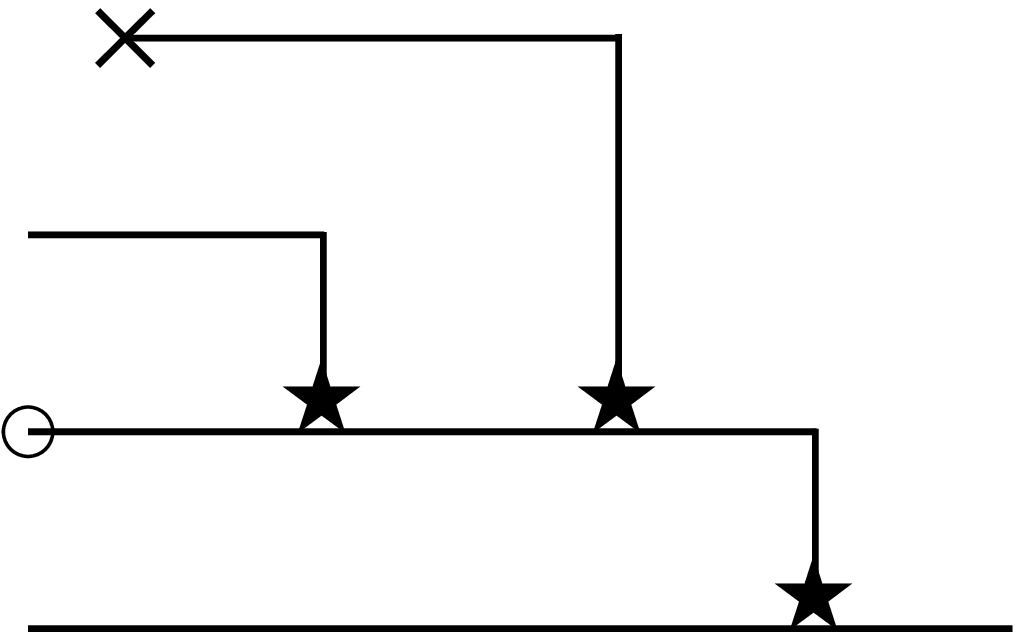}
\end{center}
\caption{\label{killed_ASG}
 Realisations of the killed ASG with viability selection started with a single individual, in the diffusion limit (first and second diagram) and the deterministic limit (third and fourth diagram).  The graph evolves in backward time, which runs from right to left. In the first and third diagram, all ancestral lineages eventually encounter mutations to type 1 (crosses), which means that the line-counting process is absorbed in 0, so the sampled individual is of type 1. In the second and fourth diagram, an ancestral lineage encounters a mutation to type 0 (a circle), which means that the line-counting process is absorbed in the cemetery state $\Delta$, so the sampled individual is of type 0.  In the absence of type-0 mutations,  the number of lines may  converge to $\infty$ in the deterministic setting.
 }
 \end{figure}

\paragraph{Diffusion limit.}
 We are guided by two elementary, but crucial insights. First, the type of an individual at present is determined by the most recent mutation along its ancestral line. Therefore, once one encounters, working back into the past, a mutation on a line, this line need not be considered any further -- it may be \emph{pruned}.  Second, type-0 individuals have priority at every branching event; the most recent beneficial mutation on a line that is still alive therefore  decides that there is at least one type-0 individual in the sample. This leads us to the following definition (see Fig.~\ref{killed_ASG}).
 
\begin{definition} \label{DefkilledASG} (killed ASG, diffusion limit).
The killed ASG in the diffusion limit starts with one line emerging from each of the   $n$ individuals in the sample. Every line branches at rate $\sigma$, every ordered pair of lines coalesces at rate 1; every line is pruned at rate $\vartheta \nu_1$. Furthermore, at rate $\vartheta \nu_0$ per line, the process is killed, that is,  reaches what we call the cemetery state $\Delta$. 
\end{definition}

This killed ASG is related to the \emph{coalescent with killing} \cite[Ch.~1.3.1]{Durrett_08}, which, in the neutral case, determines \emph{all} types in a sample; the coalescent with killing, in turn,  is \emph{Hoppe's urn model} \cite{Hoppe_84} run backward in time. Note that, in our definition, we distinguish  between \emph{pruning} (individual lines) and \emph{killing} (the entire process). Note also that the killed ASG does not yield the full type configuration of a sample, but only informs us whether or not \emph{all} individuals in the sample are of type~1.

\begin{figure}
\begin{center}
\includegraphics[ width=0.25\textwidth]{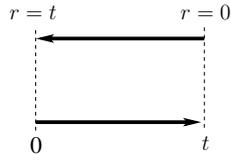}
\caption{\label{times} Forward time $t$ and backward time $r$.}
\end{center}
\end{figure}

Let now $R:= (R_r)_{r \geqslant 0}$ be the line-counting process of the killed ASG (we use the variables $t$ and $r$ throughout for forward and backward time, respectively, so $r=t$ in backward time corresponds to $t=0$ in forward time, see Fig.~\ref{times}). 
 It is a con\-tin\-uous-time Markov chain on $\NN_{\geqslant 0}\cup \{\Delta\}$. From Definition~\ref{DefkilledASG}, the process has transition rates
\begin{equation}\label{rates_R}
q^{}_R(n,n+1)  =  n \sigma, \quad
q^{}_R(n,n-1)  =  \frac{n (n-1)}{2} + n \vartheta \nu_1, \quad
q^{}_R(n,\Delta)  = n  \vartheta \nu_0
\end{equation}
for $n \in \NN_{\geqslant 0}$. The states 0 and $\Delta$ are  absorbing; all other states are transient. Note that  with probability 1 convergence to $\infty$ does not occur, since the linear birth rates are overcompensated by the quadratic death rates when the number of lines is large.  Absorption in 0  (absorption in $\Delta$) implies that (not) all individuals in the sample are of type 1. We have the following \emph{duality relation} between $X$ and $R$.
\begin{proposition}\label{duality_X_R}
Let $X$ be the Wright-Fisher diffusion and $R$ the line-counting process of the killed ASG in the diffusion limit. Defining $(1-x)^\Delta:= 0$, we then have
\begin{equation}\label{E_X_R}
\EE \big ((1-X_t)^n \mid X_0=x\big ) = \EE\big ((1-x)^{R_t} \mid R_0=n\big )
\end{equation}
for any $t \geqslant 0$, $n \in \NN_{\geqslant 0}$, and $x \in [0,1]$. 
\end{proposition}
Before proving \eqref{E_X_R} let us remark that dualities like this one are key to understanding the stochastic processes of population genetics. They come from the area of interacting particle systems (see \cite[Chaps.~3 and 4]{Liggett_10}) and often involve relations between processes forward and backward in time; see also \cite{Jansen_Kurt_14}. Proposition~\ref{duality_X_R} (which generalises Eq.~(1.5) in \cite{Athreya_Swart_05}) may be interpreted as follows. Consider  a population that has started with $X_0=x$ and is in state $X_t$ at time $t$. Take a sample of $n$ individuals from this population at time $t$.  
The  left-hand side of \eqref{E_X_R}  is the probability that all $n$ individuals are of type $1$. The statement then says that this probability may be determined by starting the killed ASG with $n$ lines at time $r=0$, running it until time $r=t$, and assigning types to each of the $R_t$ lines according to the weights $(x,1-x)$, in an i.i.d.\ fashion. All $n$ individuals in the sample are then of type 1 if and only if all $R_t$ individuals are of type 1. This is clear for the limiting case $\vartheta =0$ because then, as explained in Sec~\ref{sec:ASG}, $R_t$ when started in $n$  just counts the number of potential ancestors at (forward) time $0$ of  the sample taken at time $t$. 

Now throw mutations on the ASG. If the first mutation encountered along a given ancestral line, when proceeding back into the past, is a  type-1 mutation, then this mutation passes type 1 on to its decendants within the sample, and so this lineage need not be pursued futher back into the past.  This results in a  pruning of  lineages in the ASG at rate $\nu_1\vartheta$ per line. If on the remaining part of the ASG, when proceeding into the past, one encounters a  type-0 mutation, then this mutation passes on type 0 to its decendants within the sample, and therefore such a realisation does not count for the event that all  individuals in the sample are of type 1. This is  incorporated by a ``killing'' of the entire ASG at rate $\nu_0\vartheta$ per lineage. 

We now complement the just-stated ``graphical proof'' of the relation \eqref{E_X_R} by a more formal generator argument. For this we note that \eqref{E_X_R} can be re-expressed in terms of the semigroups $P^X$ and  $P^R$ of the processes $X$ and $R$, respectively, along with the duality function $H(x,n):= (1-x)^n$, as
\begin{align}\label{dualP}
P_t^XH(.,n)(x)  =  P_t^RH(x,.)(n).
\end{align}
Now a well-known result (Thm.~3.42 in \cite{Liggett_10}, see also \cite[Prop.~1.2]{Jansen_Kurt_14}) says that, in order to verify the duality relation \eqref{dualP}, it is enough to check the corresponding relation for the generators
\begin{align}\label{dualG}
G^XH(.,n)(x)  =  G^RH(x,.)(n).
\end{align}
In our case the two generators have the form
\begin{align*}G^Xf(x) =& \left(\frac 1 2x(1-x) \frac{\dd^2}{\dd x^2}+ \big(\sigma x(1-x)  -  \vartheta \nu_1 x+ \vartheta \nu_0(1-x)\big) \frac{\dd}{\dd x}\right)f(x),\\ 
G^Rg(n) = & \, {n\choose 2} \big (g(n-1)-g(n) \big )+\sigma n \big (g(n+1)-g(n) \big )\\&+ \vartheta \nu_1 n \big (g(n-1)-g(n) \big ) -  \vartheta \nu_0 n g(n),
\end{align*}
where the expression for $G^X$ follows from \eqref{sde} and that for $G^R$ from  \eqref{rates_R}.
The relation \eqref{dualG} now follows by a straightforward calculation. In fact, it holds true pairwise  for each of the four parts of the generators (neutral reproduction--coalescence, selection--branching, type-1 mutation--pruning, type-0 mutation--killing). For example, for the first pair (neutral reproduction--coalescence), it reads 
$$\frac 1 2x(1-x) \frac{d^2}{dx^2}(1-x)^n = {n\choose 2} \big ( (1-x)^{n-1}-(1-x)^{n} \big ),$$
which is the generator formulation of the moment duality between the classical Wright-Fisher diffusion and the line-counting process of Kingman's coalescent.  
\hfill \qed

\smallskip

Having thus completed the proof of Proposition \ref{duality_X_R}, we
let $t\to \infty$ in \eqref{E_X_R}. Assuming  $\vartheta > 0$, we observe that $R$ is eventually absorbed either in $0$ or in~$\Delta$; see Fig.~\ref{killed_ASG}. (The event $\{R_\infty=\Delta\}$ does not count for the expectation $\mathbb E((1-x)^{R_\infty})$, in line with our setting $(1-x)^\Delta := 0$ of Proposition~\ref{duality_X_R}.) Hence we see that, irrespective of the initial condition $x$, the process $X_t$ converges in distribution to an equilibrium state $X_\infty$ with
\begin{align}\label{wrightmom}
\mathbb E \big ( (1-X_\infty)^n \big ) = \mathbb P(R_\infty=0 \mid  R_0=n) =: b(n).
\end{align}
A first-step decomposition (that is, applying total probability plus the Markov property at the time of the first jump of $R$) immediately shows that the absorption probabilities $b(n)$ satisfy the recursion
\begin{align}\label{recu}
  b(n) = \frac{2 \sigma}{n-1+2\sigma + 2\vartheta} b(n+1) + \frac{n-1+2\vartheta \nu_1}{n-1+2\sigma + 2\vartheta} b(n-1),
\end{align}
complemented by the boundary conditions
\[
b(0)=1 \quad \text{and }  \; \lim_{n\to \infty} b(n)=0.
\]
The recursion \eqref{recu} is an example of a \emph{sampling recursion} as already obtained in \cite{Krone_Neuhauser_97} via a different route;
here we add the genealogical interpretation via the killed ASG (for the case $\vartheta>0$ and $0<\nu_0 < 1$). 
The limiting state $X_\infty$ has density \eqref{pi}; consequently, the $b(n)$, which by  \eqref{wrightmom} are ``moments'' of the Wright distribution, are of the form
\begin{align}\label{bnanal}
\frac{b(n)}{C}=   \int_0^1x^{2\vartheta\nu_0-1}(1-x)^{n+2\vartheta\nu_1-1} \ee^{2\sigma x}\dd x, \quad n=1,2,\ldots
\end{align}

Let us now consider $\EE(X_\infty)$, the expected proportion of type-0 individuals at stationarity, for $\nu_0 \ll 1$, as a function of the mutation rate, and in dependence of $N$. That is, we take
the approximating diffusions defined by the choice $\sigma = N s$, $\vartheta = N u$, for given $s$ and $u$. Fig.~\ref{stat_forward} (right panel) depicts the corresponding curves and illustrates their convergence to $z_\infty$, the corresponding stationary frequency in the deterministic limit according to \eqref{z_infty} as $N \to \infty$. It is now time to return to the deterministic limit and consider the killed ASG in this setting.

\paragraph{Deterministic limit.} 
In the deterministic limit, the argument is similar, but now there are no coalescence events (see Fig.~\ref{killed_ASG}, third and fourth diagram), and, for $\nu_0=0$,  the killed ASG may grow to infinite size. In what follows, we mainly rely on \cite{Baake_Cordero_Hummel}, where  details, proofs, and further results may be found.

\begin{definition} (killed ASG, deterministic limit)
The killed ASG in the deterministic limit starts with one line emerging from each of the   $n$ individuals in the sample. Every line branches at rate $s$; every line is pruned at rate $u \nu_1$; the process is killed at rate $u \nu_0$ per line. 
\end{definition}
Let us note that, due to the absence of coalescence events, this killed ASG, starting from $n$ individuals, actually consists of $n$ independent killed ASGs, each starting with a single line.

Let $R:= (R_r)_{r \geqslant 0}$ be the line-counting process of the killed ASG in the deterministic limit. It is a continuous-time Markov chain on $\NN_{\geqslant 0} \cup \{\Delta\}$,   with transition rates
\[
q^{}_{R}(n,n+1)  =  n s, \quad 
q^{}_{R}(n,n-1)  =   n u \nu_1, \quad
q^{}_{R} (n,\Delta)  = n  u \nu_0
\]
for $n \in \NN$.  The states 0 and $\Delta$ are  absorbing; all other states are transient. Absorption in 0  implies that all individuals in the sample are of type 1; absorption in $\Delta$ entails that at least one individual is of type 0.  The latter also holds if $R$ converges to $\infty$, provided $x>0$.  
The results analogous to those in the diffusion limit now read as follows.
\begin{proposition}\label{duality_det} Let $z(t; x)$ be the solution of the deterministic
mutation-selection equation \eqref{det_muse} with initial value $x\in (0,1]$, and $R$ the line-counting
process of the killed ASG in the deterministic limit. We then have
\[
\big (1-z(t;x) \big )^n = \EE\big ( (1-x)^{R_t} \mid R_0=n \big )  
\]
for  $n \in \NN_{\geqslant 0} \cup \{\Delta\}$ and $t \geqslant 0$, where  $(1-x)^\Delta=0$.
For a sample of size 1, one obtains the asymptotic behaviour
\[
  \PP \big ( \lim_{r \to \infty} R_r \in   \{\Delta,\infty\} \mid R_0=1 \big ) =  z^{}_\infty
\]
with $z^{}_\infty$ from \eqref{z_infty}.
\end{proposition}
 
Proposition~\ref{duality_det} provides an illuminating connection between the solution of the deterministic
mutation-selection equation forward in time and the stochastic killed ASG backward 
in time; indeed, it gives  a \emph{stochastic representation} of the deterministic solution. A proof is given in \cite{Baake_Cordero_Hummel}. The graphical explanation of the first statement is analogous to the diffusion limit. Let us only provide an illustrative argument for the second statement here. Let $w := \PP(R \text{ absorbs  in } 0 \mid R_0=1)$.  A decomposition according to the first step gives 
\begin{equation}\label{w}
w = \frac{u \nu_1}{u+s} + \frac{s}{u+s} w^2,
\end{equation}
where we have used that $\PP(R \text{ absorbs  in } 0 \mid R_0=2)=w^2$
due to the conditional independence of the two individuals after the branching event.
One is therefore left with a quadratic equation; its unique solution in $[0,1]$ is $w= 1-z_\infty$
with $z_\infty$ from \eqref{z_infty}. In the limiting case $\nu_0=0$, where $\Delta$ cannot be accessed, 
the bifurcation at $u=s$ in \eqref{eq_spl}  marks the dichotomy between 
the two possible fates of the birth-death process: for $u \geqslant s$, it dies out almost surely, whereas for $u<s$,  it survives with  positive probability $1-u/s$ and then grows to infinite size almost surely. This is a classical result from the theory of branching processes \cite[Ch.~III.4]{Athreya_Ney_70}: Indeed, for $\nu_1=1$, \eqref{w} is the fixed point equation $w=\varphi(w)$ for the generating function $\varphi$ of the offspring distribution of a binary Galton-Watson process with probability $u/(u+s)$ for no offspring and $s/(u+s)$ for two offspring individuals. This connection sheds new light on  \eqref{eq_spl}. Namely, let us consider the killed ASG starting from a single  individual sampled from the equilibrium population (at some late time $t$, say). Then $R_0=1$, and on the event $\{R_r \to 0$  for $r \to \infty\}$  the sampled individual  is of type 1. In contrast, on the event $\{R_r \to \infty\}$  any (even the smallest) positive value of $x$ suffices to ensure that a type 0  is assigned to at least one ancestral line, which guarantees that the individual sampled from the equilibrium population is of type 0.

\begin{remark}
Due to the independence of the individuals in the killed ASG in the deterministic limit, the probability for arbitrary type configurations of a sample of size $n$ is easily determined via $n$ independent killed ASGs. This is different in the diffusion limit, where individuals are dependent via common ancestry; we therefore only ask whether or not all $n$ individuals are of type~1 in this case. In fact, the construction may be extended to yield arbitrary type configurations, but this requires additional effort.
\end{remark}

\section{The pruned lookdown ASG and the ancestral type distribution}
\label{sec:pruned}
Let us now turn to a graphical construction of the type of the \emph{ancestor} at time 0 of an individual  chosen randomly at a fixed later time $t$. This is a more involved problem than identifying the (stationary) type distribution of the forward process, because we now must identify the parental branch (incoming or continuing, depending on the type) at every branching event, which requires  nested case distinctions. In \cite{Lenz_Kluth_Baake_Wakolbinger_15}, we have overcome this problem, in the case of fecundity selection, by introducing an \emph{ordering} of the lines (analogous to the one already used in  Figs.~\ref{IC}--\ref{killed_ASG}); this was combined with a pruning procedure, which, upon mutation, eliminates lines that can never be ancestral. Moreover, we place the lines of the ordered graph on consecutive levels, starting at level 1. This bears elements of the aforementioned \emph{lookdown construction} \cite{Donnelly_Kurtz_99}, which are thus combined with the (pruned) ASG, hence the name \emph{pruned lookdown ASG}. We now describe the construction, this time starting with the deterministic limit.

\paragraph{Deterministic limit.} We restrict ourselves to a sample of size $n=1$, since ancestries are independent due to the absence of coalescence events. The construction follows \cite{Baake_Cordero_Hummel,Cordero_17} and is extended to the case with viability selection.
\begin{definition} \label{pLDASG_det} (pruned lookdown ASG, deterministic limit)
The pruned lookdown ASG in the deterministic limit starts with one line  at time $r=0$  and proceeds in direction of increasing $r$. At each time $r$,  the graph consists of a finite number $L_r$ of lines. 
The lines are numbered by the integers $1,\ldots, L_r$, to which we refer as \emph{levels}. The process then evolves via the following transitions
(see Figs.~\ref{elements_ASG_det} and \ref{example_ASG_det}).
\begin{figure}
\begin{center}
\includegraphics[width=0.5\textwidth]{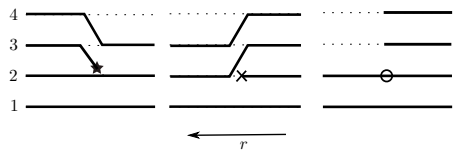} 
\end{center}
\caption{\label{elements_ASG_det} Transitions of the pruned lookdown ASG with viability selection in the deterministic limit: branching, pruning, and killing. 
}
\end{figure}

\begin{figure}
\begin{center}
\includegraphics[width=0.6\textwidth]{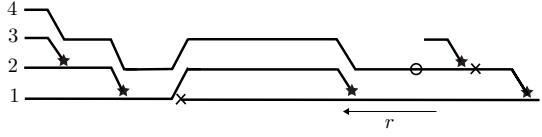}
\end{center}
\caption{\label{example_ASG_det} A cut-out of a realisation  of the pruned lookdown ASG with viability selection in the deterministic limit.}
\end{figure}

\begin{enumerate} 
	\item
	Every line $i \leqslant L_r$ branches at rate $s$ and  a new line, namely the incoming branch, is inserted. In the case with fecundity selection, the insertion is at level $i$ and all lines at levels $k\geqslant i$ 
	are pushed one level upward to $k+1$;  in particular, the continuing branch is shifted from level $i$ 
	to  $i+1$. With viability selection, the continuing branch remains at level $i$,  the incoming branch is inserted  at level $i+1$, and  all lines at levels $k\geqslant i+1$ 
	are pushed one level upward to $k+1$. In any case,  $L_r$ increases to  $L_{r}+1$. 
	\item 
	Every line $i \leqslant L_r$ experiences deleterious mutations  at rate $u \nu_1$. If $i=L_r$, nothing happens. If $i<L_r$,   the line at level $i$ is pruned,  and the lines above it  slide down to `fill the gap', rendering the transition from $L_r$ to  $L_{r}-1$. 
	\item 
	Every line $i \leqslant L_r$ experiences beneficial mutations  at rate $u \nu_0$.  All the lines at levels $> i$ are pruned, resulting in a transition from $L_r$ to $i$. Thus, no pruning happens if a beneficial mutation occurs on level $L_r$.

%	 (if $i=L_{r} $,  in fact nothing happens).
\end{enumerate}
\end{definition}
Let us explain the idea behind the process. The ordering, brought about by the placement of each incoming line immediately below  or above the continuing line as anticipated in Section~\ref{sec:ASG}, entails that the levels reflect the hierarchy according to the pecking order. To see this, consider first the case without mutation and hence without pruning. It is then clear that every line has, at some  point in the forward direction of time, priority over the line above it --- unless it is  the top line, which is continuing in all branching events in which it is involved. As a consequence, there is a hierarchy from bottom to top in the sense that the level of the ancestral line at time 0 is either the lowest type-0 level at time 0 or, if all $L_0$ lines are of type 1, it is  level~$L_0$. 

Now add in the mutations. Encountering a deleterious mutation on a line at a level below $L_r$ entails that this line will not be ancestral at the branching event at which it has priority over the line above it; it therefore need not be considered as a potential ancestor any further and may be pruned. In contrast, a deleterious mutation at the top level does not lead to pruning since the top line will be ancestral regardless of its type, provided all lines below it carry type 1; the top line is therefore called \emph{immune} (to deleterious mutations). A beneficial mutation implies that no line above the one that carries the mutation can be ancestral, which results in the corresponding pruning action. Since all pruning operations preserve the order of the existing lines, the pecking order is retained throughout.

As a consequence of  Definition~\ref{pLDASG_det}, the process $L=(L_r)_{r \geqslant 0}$
has transition rates
\begin{equation}\label{rates_det_L}
q^{}_L(n,n+1) = n s, \quad q^{}_L(n,n-1) = (n-1) u \nu^{}_1 + u \nu^{}_0 \one\{n>1\}, \quad q^{}_L(n,n-\ell) = u \nu^{}_0, 
\end{equation}
$2 \leqslant \ell <  n, \quad n \in \NN$. 
Let us summarise its asymptotic behaviour (following \cite{Baake_Cordero_Hummel}).
\begin{proposition}\label{L_infty}
For the pruned lookdown ASG in the deterministic limit, we have
\begin{enumerate}
\item 
For $s=0$, one has $L_r \equiv 1$, so, in particular, $L_\infty := \lim_{r \to \infty} L_r  = 1$. 
\item For $s>0$, $u \leqslant s$, and $\nu_0=0$, $L_\infty=\infty$ (almost surely for $u<s$, in probability for $u=s$).
\item
For $s>0$ and either $\nu_0>0$ or $u>s$, the process  attains a  stationary distribution; the corresponding random variable, denoted again by $L_\infty$, has the geometric distribution $\Geo(1-p)$ with parameter 
\begin{equation}\label{geo_p}
p = \begin{cases} \frac{1}{2} \Big (\frac{u+s}{u \nu_1} - \sqrt{\Big ( \frac{u+s}{u \nu_1})^2 - 4 \frac{s}{u \nu_1}} \Big ), & \nu^{}_1 > 0, \\
\frac{s}{u+s}, & \nu^{}_1 = 0. \end{cases}
% \qedhere
\end{equation}
\end{enumerate}
\end{proposition}
See   \cite{Baake_Cordero_Hummel} for a proof that relies on the graphical construction and also provides insight into the property of `no memory'  that leads to the geometric distribution.
Note that, for $\nu_1>0$, one has $p= \frac{s}{u \nu_1} (1-z_\infty)$ with  $z_\infty$ of  \eqref{z_infty}. Note also that the distribution of $L_\infty$ in the cases $s=0$ as well as $s>0, u \leqslant s, \nu_0=0$ may be seen as degenerate cases of $\text{Geo}(1-p)$. Namely, for $s=0$, one has $p=0$ (in agreement with \eqref{geo_p}), which means immediate success, in line with $L_\infty=1$. In contrast, for $s>0, u \leqslant s, \nu_0=0$, we set $p=1$, which is consistent with an infinite number of trials. 

\begin{figure}
\begin{center}
\includegraphics[ width=0.48\textwidth]{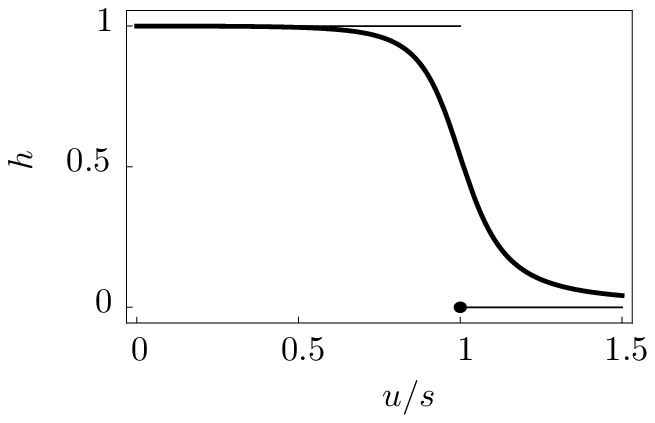} \quad
\includegraphics[ width=0.48\textwidth]{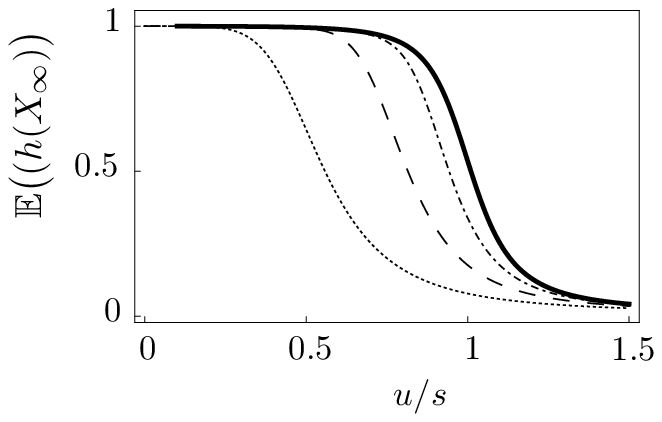} 
\end{center}
\caption{
\label{stat_backward}
The stationary proportion of a type-0 ancestor in the deterministic limit (left) and the corresponding expectation in the diffusion approximation (right) for $s=0.001$ as a function of the mutation rate. Parameters as in Fig.~\ref{stat_forward}. The bold line is $h(z_\infty)$ for $\nu_0=0.005$ in both panels.
}
\end{figure}

Consider now the sampling of the potential ancestors' types at time 0. Due to the pecking order, the true ancestor at time 0 of an individual at time $t$ is of type 1 if and only if all  $L_t$ potential ancestors of the individual are assigned type 1 when sampled from the distribution with weights $(x, 1-x)$ in an i.i.d.\ fashion. We are particularly interested in the limit $t \to \infty$, that is, in the type of the ancestor of a random individual sampled from the equilibrium distribution, given that the initial frequency of the beneficial type was $x$.
The probability $h(x)$ that this ancestor is of type 0 is then given by the probability of at least one success in a random number of $L_\infty$ coin tosses, each with success probability $x$,
as summarised in the following theorem, once more from \cite{Baake_Cordero_Hummel,Cordero_17}.
\begin{theorem}\label{thm_anc_det}
Let $J_t$ be the type of the ancestor at time $0$ of an individual randomly sampled from the population at time $t$ in the deterministic limit. For $x \in [0,1]$, we then have
\[
 h(x) := \lim_{t \to \infty} \PP(J_t=0 \mid X_0=x) = \sum_{n \geqslant 0} x (1-x)^n a^{}_n,
\]
where $a^{}_n := \PP(L_\infty > n) = p^n$ with $p$ from Proposition~\ref{L_infty} (including the limiting cases). 
\end{theorem}
Let us explain in words what Theorem~\ref{thm_anc_det} tells us. In the neutral case (that is,  $s=0$), we have $p=0$ and hence $a_0=1$ and $a_n=0$ for all $n \geqslant 1$. Hence $h(x)=x$, so there is no bias towards one of the two types. In contrast, for $s>0$, we have $p>0$ (increasing in $s$) and so $a_n>0$ (also increasing in $s$) for all  $n \geqslant 0$. This explains that, and how, selection introduces a bias towards the beneficial type in the ancestry: It increases the number of potential ancestors, thus providing more chances  for the ancestor to be of type $0$.  
It is particularly interesting to start from a stationary population, that is, $x=z_\infty$ of \eqref{eq_spl}. Fig.~\ref{stat_backward} (left panel) shows $h(z_\infty)$ as a function of $u$; comparing this with the left panel of Fig.~\ref{stat_forward} illustrates the bias in an impressive way. Consider, in particular,
 $s>0$ and the limiting case of $\nu_0=0$.  Then
Proposition \ref{L_infty}  tells us that
 $L_\infty   = \infty$ for $u \leqslant s$, whereas $L_\infty$ follows $\text{Geo}(s/u)$ if $u>s$. 
 With the coin-tossing interpretation of $h(x)$ given before Theorem \ref{thm_anc_det}, we see immediately that $h(z_\infty) = 1$ for $u < s$ (since then $z_\infty  >0$ and with probability 1 there occurs a success in an infinite number of trials), and that $h(z_\infty) = 0$ for $u \geqslant s$ (since then $z_\infty = 0$ and with probability 1 there is no success even in an infinite number of trials).

Still considering $\nu_0=0$, we therefore see that the error threshold for the stationary type distribution  \eqref{eq_spl} is accompanied by a more drastic effect at the level of the ancestral type distribution, which \emph{jumps} from a point measure on 0 to a point measure on 1. This behaviour was found earlier  \cite{Hermisson_etal_02} via analysis of the first-moment generator of the multitype branching process mentioned the introduction, but now appears in a new light. Namely, with an infinite number of Bernoulli trials (for $u < s$), any positive proportion of beneficial individuals at time $0$ will  guarantee that the ancestor at time 0 of a randomly chosen individual from the limiting distribution at $t \to \infty$ is of type 0. In contrast, for $u \geqslant s$, one samples finitely or infinitely many potential ancestors from a pure type-1 population; this yields a type-1 ancestor with probability 1.

\begin{figure}
\begin{center}
\includegraphics[width=0.5\textwidth]{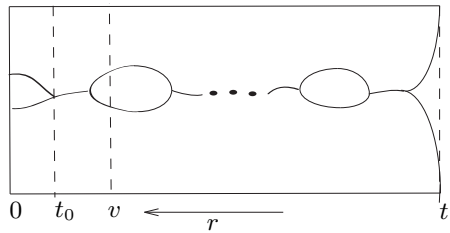}
\end{center}
\caption{\label{ancline}
Schematic sketch of a realisation of an equilibrium ASG (evolving in time $r$), with $t_0$ denoting its smallest nonnnegative (forward) time of a bottleneck.
%Ancestral lines in a population at time $t$, when traced back into the past, coalesce into a single line of descent. The bold lines symbolise the true ancestral lines emerging from distinct individuals at time $t$, whereas the thin black ‘envelope’ hints at the ancestry of the entire population.
}
\end{figure}

\paragraph{Diffusion limit.} 
In the diffusion limit, ancestral lines can also coalesce.  As with the killed ASG in the diffusion limit, this entails that individuals  no longer have independent ancestries. 
Rather, for large enough $t$, the potential ancestral lines  of an infinite sample of individuals drawn from the population at time $t$ will, on their way back to time $0$, eventually coalesce into a single line (and then branch again). In other words: the ASG on its way back into the past has \emph{bottlenecks}, that is, with probability 1 it repeatedly returns to a state in which it consists of a single line.  Let $t_0$ be the smallest among all the non-negative (random) times at which there is a bottleneck of the  ASG, see Fig.~\ref{ancline} for an illustration.  Then, for determining the type of that individual at time $0$ that is ancestral to the entire population at the (late) time $t$, it suffices to consider the ASG between (forward) times $0$ and $t_0$. 

For $t \to \infty$, the restriction of the ASG  to any (forward) time interval $[0,\upsilon]$ will stabilise in distribution, rendering in the limit the so-called {\em equilibrium ASG}. (In the limit $t\to \infty$ it plays no role whether the ASG is started from infinitely many lines or, say, from a single line at time $t$.)

\begin{figure}
\begin{center}
\includegraphics[width=0.9\textwidth]{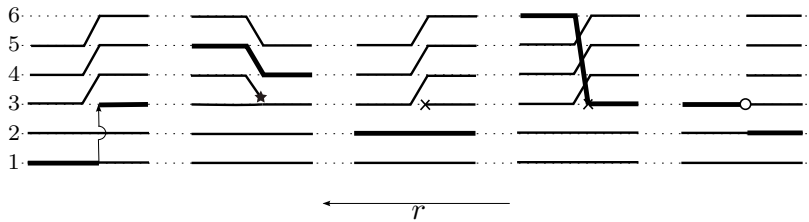} 
\end{center}
\caption{\label{elements_diff} Transitions of the pruned lookdown ASG with viability selection in the diffusion limit.
The immune line is marked in bold. From left to right: coalescence; branching; pruning due to a deleterious mutation outside the immune line; relocation due to a deleterious mutation on the immune line; pruning due to a beneficial mutation.
}
\end{figure}

\begin{figure}
\begin{center}
\includegraphics[width=0.6\textwidth]{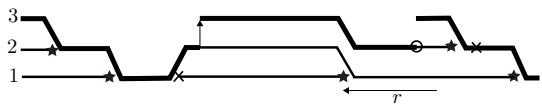}
\end{center}
\caption{\label{realisation_diff} A cut-out of a realisation  of the pruned lookdown ASG with viability selection in the diffusion limit. The immune line is marked in bold.}
\end{figure}

In order to compute the distribution of the type of the common ancestor at time $0$, we may therefore consider the equilibrium ASG and, after introducing mutations along the ASG, its ordered and pruned version, the pruned lookdown ASG in equilibrium. As in the deterministic limit, the pruned lookdown ASG has $L_r$ lines at backward time $r$, numbered by their levels. As will be detailed below, branching  events are as in the deterministic limit, with $s$ replaced by $\sigma$. In addition, there are coalescence events at rate $1$ per  pair of lines. Pruning (with rate $u$ replaced by $\vartheta$) is similar to before, with one important modification. At any given time, there is again exactly one immune line that is unaffected by deleterious mutations (and is ancestral if all lines are of type 1). But this line need no longer be the top line (this is due to the coalescence events, which can move the line downwards). The precise rules where this line is located and how it is relocated by the various events are derived in \cite{Lenz_Kluth_Baake_Wakolbinger_15} and \cite{Baake_Lenz_Wakolbinger_16}. The result is summarised, and extended to viability selection, in the following definiton (see Figs.~\ref{elements_diff} and \ref{realisation_diff}).

\begin{definition} (pruned lookdown ASG, diffusion limit)
The pruned lookdown ASG in the diffusion limit starts with one line  at time $r=0$.   At each time $r>0$,  the graph consists of a finite number $L_r$ of lines, one of which is distinguished (and is called immune). 
The lines are numbered by the levels $1,\ldots, L_r$; the level of the immune line is $M_r \leqslant L_r$. The process  evolves via the following transitions.
\begin{enumerate} 
	\item
	Every line $i \leqslant L_r$ branches at rate $\sigma$, and then  a new line, namely the incoming branch, is inserted. In the case with fecundity selection, the insertion is at level $i$ and all lines at levels $k\geqslant i$ 
	are pushed one level upward to $k+1$;  in particular, the continuing branch is shifted from level $i$ 
	to  $i+1$. With viability selection, the continuing branch remains at level $i$,  the incoming branch is inserted  at level $i+1$, and  all lines at levels $k\geqslant i+1$ 
	are pushed one level upward to $k+1$. In any case,  $L_r$ increases to  $L_{r}+1$. If $M_r \geqslant i$, then $M_r$ increases to  $M_{r}+1$; otherwise, it remains unchanged. 
        \item Every ordered pair of lines $(i,j)$, $i<j \leqslant L_r$, coalesces at rate 1. The remaining lines are relocated to `fill the gap' while retaining their original order; thus $L_r$ decreases by one. The immune line follows the line on level $M_r$.
%that is, it decreases to  $M_{r}-1$ if $M_r > i$; otherwise, it remains unchanged.
	\item 
	Every line $i \leqslant L_r$ experiences deleterious mutations  at rate $u \nu_1$.   If  $i \neq M_r$, then the line at level $i$ is pruned,  and the remaining lines (including the immune line) are relocated to `fill the gap' (again in an order-preserving way), rendering the transition of $L_r$  to $L_r-1$. If, however, $i = M_r$, then the line affected by the mutation is not pruned but relocated to the currently highest level, that is, $M_r$ increases to $L_r$. All  lines above  $i$ are shifted one level down, so that the gaps are filled, and in this case $L_r$ remains unchanged.
	\item 
	Every line $i \leqslant L_r$ experiences beneficial mutations  at rate $u \nu_0$.  All  lines at levels $> i$ are pruned, resulting in a transition from $L_r$ to $i$. The immune line is relocated to level $i$.
\end{enumerate}
\end{definition}

The transition rates for $L$ follow directly from this definition and read
\begin{equation}\label{q_L}
\begin{split}
& q^{}_L(n,n+1) = n s, \quad q^{}_L(n,n-1) = \frac{1}{2} n (n-1) + (n-1) u \nu^{}_1 + u \nu^{}_0 \one\{n>1\},  \\
& q^{}_L(n,n-\ell) = u \nu^{}_0, \quad 2 \leqslant \ell \leqslant  n, \quad n \in \NN;
\end{split}
\end{equation}
apart from the different parameter scaling, these rates differ from those in the deterministic limit   \eqref{rates_det_L}  only via the additional   coalescence events. Again, the process has a unique stationary distribution, which corresponds to $L_\infty$, and we  have a result analogous to Theorem~\ref{thm_anc_det}.
\begin{theorem}
Let $J_t$ be the type of the ancestor at time $0$ of a random individual at time $t$ in the diffusion limit. For $x \in [0,1]$, we then have
\[
h(x):= \lim_{t \to \infty} \mathbb{P}(J_t=0\mid X_0=x) =\sum_{n\geqslant 0} x(1-x)^{n}a_n,
\]
	where the $a_n= \PP(L_\infty>n)$ are the unique solution to \emph{Fearnhead's recursion}
\begin{equation}\label{fearnhead}
 \Big [ \frac{1}{2} (n+1)+\sigma+\vartheta\Big ]  a^{}_n = \Big [ \frac{1}{2} (n+1)+\vartheta \nu^{}_1\Big ]  a^{}_{n+1} + \sigma a^{}_{n-1}, \quad n\geqslant 1,
\end{equation}
with $a_0=1$ and $\lim_{n \to \infty} a_n=0$.
\end{theorem}
 Predecessors of this result go back to Fearnhead  \cite{Fearnhead_02} and Taylor \cite{Taylor_07}; probabilistic proofs were given in \cite{Lenz_Kluth_Baake_Wakolbinger_15,Baake_Lenz_Wakolbinger_16}. Here, we only reprove \eqref{fearnhead}, in a way that is simpler and more elegant than previous versions, and is directly based on the tail probabilities and the graphical construction. The following proof is a straightforward extension of the  proof for the deterministic limit (see Prop.~11 and Fig.~6 in \cite{Baake_Cordero_Hummel}).

%A proof, based on the stationarity condition forward in time based on the rates \eqref{q_L}, was given in \cite{Lenz_Kluth_Baake_Wakolbinger_15}. 
\begin{proof}[of Fearnhead's recursion]
 Fix $r>0$.  Let $T_*, T_\uparrow, T_\circ$ and $T_\times$ be the times of the (in the direction of $L$) most recent selective, coalescence, beneficial, and deleterious mutation event \emph{on the first $n$ levels}  (note that the $T$'s depend on $r$). Set $T:=\max\{T_*, T_\uparrow,T_\circ,T_\times\}$. Then 
\[
\begin{split}
\PP(L_r>n) & =\PP(L_r>n, T=T_*)+\PP(L_r>n, T=T_\uparrow) \\
          & +\PP(L_r>n, T=T_\times)+\PP(L_r>n, T=T_\circ).
\end{split}
\]
Reading each transition in Fig.~\ref{elements_ASG_det} from left to right, one concludes the following. If $T=T_*$, then $L_r>n$ if and only if $L_{T-}>n-1$; here, $L_{T-} := \lim_{w\nearrow T} L_w$, that is, the state `just before' the jump.  If $T=T_\uparrow$, then $L_r>n$ if and only if $L_{T-}>n+1$. If $T=T_\times$, then $L_r>n$ if and only if $L_{T-}>n+1$. Note that the latter also holds if the line that mutates is immune (in which case it is not pruned, but relocated to the top level). The case $T=T_\circ$  contradicts $L_r>n$, so $\PP(L_r>n, T=T_\circ)=0$.  Hence, on $\{L_r > n\}$, the probabilities for the most recent event to be a branching, a coalescence, or a deleterious mutation are $\sigma/( (n+1)/2 +\sigma+\vartheta)$, $(1/2) (n+1) /( (n+1)/2 +\sigma+\vartheta)$, and $\vartheta \nu_1 / ( (n+1)/2 +\sigma+\vartheta)$, respectively, and so 
\[
\begin{split}
\Big [ \frac{1}{2} (n+1)&+\sigma+\vartheta \Big ]  \PP(L_r>n)  = \sigma \PP(L_{T-}>n-1 \mid T=T_*) \\
& + \frac{1}{2} (n+1) \PP(L_{T-}>n+1 \mid T=T_*)  
 +\vartheta \nu_1 \PP(L_{T-}>n+1 \mid T=T_\times).
\end{split}
\]
But $L_{T-}$ is independent of what happens at time $T$, since this is in the future (in $r$-time). The claim thus follows by taking $r\to \infty$ on both sides. \qed
\end{proof}

As in the deterministic limit, let us finally consider, instead of a given type frequency $x$ at time $0$, the equilibrium type frequency, which here is the random variable $X_\infty$ of Section~\ref{sec:killed}. This amounts to starting in a stationary type distribution at time $0$. The resulting expected ancestral frequency of the beneficial types, $\mathbb{E}(h(X_\infty))$, is illustrated in the right panel of Fig.~\ref{stat_backward}, again for various population sizes in the corresponding diffusion approximation.

%\subsection{Subsection title}
%\label{sec:2}
%as required. Don't forget to give each section
%and subsection a unique label (see Sect.~\ref{sec:1}).
%\paragraph{Paragraph headings} Use paragraph headings as needed.
%\begin{equation}
%\end{equation}

%\begin{figure}
%  \includegraphics{example.eps}
%\end{figure}
%\begin{figure*}
%  \includegraphics[width=0.75\textwidth]{example.eps}
%\end{figure*}
%\begin{table}
%\begin{tabular}{lll}
%\hline\noalign{\smallskip}
%\end{tabular}
%\end{table}

\begin{acknowledgements}
It is our pleasure to thank Fernando Cordero, Sebastian Hummel, and Ute Lenz for fruitful discussions. This project received financial support from Deutsche Forschungsgemeinschaft (Priority Programme SPP 1590 \emph{Probabilistic Structures in Evolution}, grants no. BA 2469/5-1 and WA 967/4-1).
\end{acknowledgements}



\end{document}